\documentclass[11pt,letterpaper]{article}
\usepackage{amsthm}
\usepackage{amsmath}
\usepackage{amsfonts}
\usepackage{latexsym}
\usepackage{geometry}
\usepackage{enumerate}
\usepackage{csquotes}

\emergencystretch=\maxdimen
\title{A note on Perelman's no shrinking breather theorem}
\author{Yongjia Zhang}
\begin{document}
\maketitle

As an application of his entropy formula, Perelman \cite{perelman2002entropy} proved that every compact shrinking breather is a shrinking gradient Ricci soliton. We give a proof for the complete noncompact case by using Perelman's $\mathcal{L}$-geometry. Our proof follows the argument in Lu and Zheng \cite{lu2017new} of constructing an ancient solution, and removes a technical assumption made by them.  For other proofs with additional assumptions, please refer to Zhang \cite{zhang2014no}.
\\

\section{Introduction}

After showing that the Ricci flow is the gradient flow of the $\mathcal{F}$ functional
\begin{eqnarray*}
\mathcal{F}(g,f):=\int_M\left(|\nabla f|^2+R\right)e^{-f}dg,
\end{eqnarray*}
Perelman \cite{perelman2002entropy} indicated that the Ricci flow on a manifold $M$ can be regarded as an orbit in the space
\begin{eqnarray*}
\text{Met}(M)\Big/\text{Diff},
\end{eqnarray*}
where $\text{Met}(M)$ stands for the space of all the Riemannian metrics on $M$ and $\text{Diff}$ represents all the self-diffeomorphisms on $M$. The breathers are the periodic orbits in this space.

\newtheorem{Definition_breather}{Definition}
\begin{Definition_breather}
A metric $g(t)$ evolving by the Ricci flow on a Riemannian manifold $M$ is called a breather, if for some $t_1<t_2$, there exists an $\alpha>0$, a diffeomorphism $\phi:M\rightarrow M$, such that $\alpha g(t_1)=\phi^*g(t_2)$. If $\alpha =1$, $\alpha <1$, or $\alpha >1$, then the breather is called steady, shrinking, or expanding, respectively.
\end{Definition_breather}

As a special case of the periodic orbits, the Ricci solitons, moving by diffeomorphisms, are the static orbits in the space $\text{Met}(M)/\text{Diff}$.

\newtheorem{Definition_soliton}[Definition_breather]{Definition}
\begin{Definition_soliton}
A gradient Ricci soliton is a tuple $(M,g,f)$, where $(M,g)$ is a Riemannian manifold and $f$ is a smooth function on $M$ called the potential function, satisfying
\begin{eqnarray*}
Ric+\nabla^2f=\frac{\lambda}{2}g,
\end{eqnarray*}
where $\lambda=0$, $\lambda=1$, or $\lambda=-1$, corresponding to the cases of steady, shrinking, or expanding solitons, respectively.
\end{Definition_soliton}

It is well understood that when moving by the $1$-parameter family of diffeomorphisms generated by the potential function, along with a scaling factor, the pull-back metric on the soliton satisfies the Ricci flow equation, and this Ricci flow is called the \emph{canonical form} of the Ricci soliton; one may refer to \cite{chow2007ricci} for more details.
\\

Perelman proved that on a closed manifold, any periodic orbit in $\text{Met}(M)/\text{Diff}$ must be static.

\newtheorem{Perelman_nobreather}[Definition_breather]{Theorem}
\begin{Perelman_nobreather}[Perelman's no breather theorem]
A steady, shrinking, or expanding breather on a closed manifold is (the canonical form of) a steady, shrinking, or expanding gradient Ricci soliton, respectively. In particular, in the steady or expanding case, the breather is also Einstein.
\end{Perelman_nobreather}

We extend the no shrinking breather theorem to the complete noncompact case.

\newtheorem{Theorem_Main}[Definition_breather]{Theorem}
\begin{Theorem_Main}\label{Theorem_Main}
Every complete noncompact shrinking breather with bounded curvature is (the canonical form of) a shrinking gradient Ricci soliton.
\end{Theorem_Main}

Our main technique is the $\mathcal{L}$-geometry, one of the two monotonicity formulae on the Ricci flow found by Perelman. In section 2 we give a brief introduction to the $\mathcal{L}$-functional. In section 3 we prove Theorem \ref{Theorem_Main}.
\\

\section{Perelman's $\mathcal{L}$-geometry}

The definitions and results in this section can be found in Perelman \cite{perelman2002entropy} and Naber \cite{naber2010noncompact}. We consider a \emph{backward} Ricci flow $(M,g(\tau))$, $\tau\in[0, T]$, satisfying
\begin{eqnarray}\label{backward_RF}
\frac{\partial}{\partial\tau}g(\tau)=2 Ric(g(\tau)).
\end{eqnarray}
Let $\gamma(\tau):[0,\tau_0]\rightarrow M$ be a smooth curve, then the $\mathcal{L}$-functional of $\gamma$ is defined by
\begin{eqnarray}
\mathcal{L}(\gamma):=\int_0^{\tau_0}\sqrt{\tau}\Big(R(\gamma(\tau),\tau)+|\dot{\gamma}(\tau)|^2_{g(\tau)}\Big)d\tau.
\end{eqnarray}
The \emph{reduced distance} between two space-time points $(x_0,0)$, $(x_1,\tau_1)$, where $\tau_1>0$, is defined by
\begin{eqnarray} \label{reduced_distance}
l_{(x_0,0)}(x_1,\tau_1):=\frac{1}{2\sqrt{\tau_1}}\inf_\gamma\mathcal{L}(\gamma),
\end{eqnarray}
where the $\inf$ is taken among all the (piecewise)  smooth curves $\gamma: [0,\tau_1]\rightarrow M$, such that $\gamma(0)=x_0$ and $\gamma(\tau_1)=x_1$. When regarded as a function of $(x_1,\tau_1)$, we call $l_{(x_0,0)}(\cdot,\cdot)$ the \emph{reduced distance based at} $(x_0,0)$. When the base point is understood, we also write $l_{(x_0,0)}$ as $l$. It is well known that the \emph{reduced volume} based at $(x_0,0)$
\begin{eqnarray}
\mathcal{V}_{(x_0,0)}(\tau):=\int_M(4\pi\tau)^{-\frac{n}{2}}e^{-l_{(x_0,0)}(\cdot,\tau)}dg(\tau)
\end{eqnarray}
is monotonically decreasing in $\tau$. We often write $\mathcal{V}_{(x_0,0)}(\tau)$ as $\mathcal{V}(\tau)$ for simplicity. We also remark here that the integrand $(4\pi\tau)^{-\frac{n}{2}}e^{-l}$ of the reduced volume is a subsolution to the conjugate heat equation
\begin{eqnarray*}
\frac{\partial}{\partial\tau}u-\Delta u+Ru=0,
\end{eqnarray*}
in the barrier sense or in the sense of distribution.
\\

Now we consider an ancient solution $(M,g(\tau))$, where $\tau\in[0,\infty)$ is the backward time. The Type I condition is the following curvature bound.

\newtheorem{Definition_Type_I}[Definition_breather]{Definition}
\begin{Definition_Type_I}
An ancient solution $(M,g(\tau))$, where $\tau\in[0,\infty)$ is the backward time, is called Type I if there exists $C<\infty$, such that
\begin{eqnarray*}
|Rm|(\tau)\leq\frac{C}{\tau},
\end{eqnarray*}
for every $\tau\in(0,\infty)$.
\end{Definition_Type_I}

To ensure the existence of a smooth limit, the $\kappa$-noncollapsing condition is often required.

\newtheorem{Definition_noncollapsing}[Definition_breather]{Definition}
\begin{Definition_noncollapsing}
A backward Ricci flow is called $\kappa$-noncollapsed, where $\kappa>0$, if for any space-time point $(x,\tau)$, any scale $r>0$, whenever $|Rm|\leq r^{-2}$ on $B_{g(\tau)}(x,r)\times[\tau,\tau+r^2]$, it holds that $\operatorname{Vol}_{g(\tau)}(B_{g(\tau)}(x,r))\geq\kappa r^n$.
\end{Definition_noncollapsing}

We will use the following theorem of Naber \cite{naber2010noncompact}.

\newtheorem{Theorem_naber}[Definition_breather]{Theorem}
\begin{Theorem_naber}[Asymptotic shrinker for Type I ancient solution] \label{asymoptotic_shrinker}
Let $(M,g(\tau))$, where $\tau\in[0,\infty)$ is the backward time, be a Type I $\kappa$-noncollapsed ancient solution to the Ricci flow. Fix $x_0\in M$. Let $l$ be the reduced distance based at $(x_0,0)$. Let $\{(x_i,\tau_i)\}_{i=1}^\infty\subset M\times(0,\infty)$ be such that $\tau_i\nearrow\infty$ and
\begin{eqnarray} \label{l_bound}
\sup_{i=1}^\infty l(x_i,\tau_i)<\infty.
 \end{eqnarray}
 Then $\{(M,\tau_i^{-1}g(\tau\tau_i),(x_i,1))_{\tau\in[1,2]}\}_{i=1}^\infty$ converges, after possibly passing to a subsequence, to the canonical form of a shrinking gradient Ricci soliton.
\end{Theorem_naber}

\noindent\textbf{Remark 1:} In Naber's original theorem, he fixes the base points $x_i\equiv x_0$. However, it is easy to observe from his proof that so long as (\ref{l_bound}) holds, all the estimates of $l$ also hold in the same way as in his case. Hence one may apply the blow-down shrinker part of Theorem 2.1 in \cite{naber2010noncompact} to the sequence of space-time base points $(x_i,\tau_i)$ and the scaling factors  $\tau_i^{-1}$.
\\

\noindent\textbf{Remark 2:} The estimates for $l$ and the monotonicity formula for $\mathcal{V}$ in \cite{naber2010noncompact} do not depend on the noncollapsing condition. According to Hamilton \cite{hamilton1995compactness}, if the noncollapsing assumption is replaced by
\begin{eqnarray} \label{inj}
\inf_{i=1}^\infty\operatorname{inj}_{\tau_i^{-1}g(\tau_i)}(x_i)>\delta,
\end{eqnarray}
where $\operatorname{inj}_g(x)$ stands for the injectivity radius of the metric $g$ at the point $x$, and $\delta>0$ is a constant, then the conclusion of Theorem \ref{asymoptotic_shrinker} still holds.
\\

\section{Proof of the main theorem}

Following the argument in Lu and Zheng \cite{lu2017new}, we construct a Type I ancient solution to the Ricci flow starting from a given shrinking breather. After scaling and translating in time, we consider the \emph{backward} Ricci flow $(M,g_0(\tau))_{\tau\in[0,1]}$, where $g_0(\tau)$ satisfies (\ref{backward_RF}), such that there exists $\alpha\in(0,1)$ and a diffeomorphism $\phi:M\rightarrow M$, satisfying
\begin{eqnarray}\label{breather_condition}
\alpha g_0(1)=\phi^*g_0(0).
\end{eqnarray}
Furthermore, we let $C<\infty$ be the curvature bound, that is,
\begin{eqnarray} \label{bounded_condition}
\sup_{M\times[0,1]}|Rm|(g(\tau))\leq C.
\end{eqnarray}
\\

For notational simplicity, we define
\begin{eqnarray*}
\tau_i=\sum_{j=0}^i\alpha^{-j},
\end{eqnarray*}
where $i=0,1,2,...$ Apparently, $\tau_i\nearrow\infty$ since $\alpha\in(0,1)$, and we can find a $C_0<\infty$ depending only on $\alpha$ (for instance, one may let $C_0=(1-\alpha)^{-1}$), such that
\begin{eqnarray}\label{time_bound}
\alpha^{-i}\leq\tau_i\leq C_0\alpha^{-i}, \text{ for every }  i\geq 0.
\end{eqnarray}
 For each $i\geq 1$, we define a Ricci flow
\begin{eqnarray}\label{piece_metric}
g_i(\tau):=\alpha^{-i}(\phi^{i})^*g_0\left(\alpha^i(\tau-\tau_{i-1})\right),\ \text{where }\tau\in[\tau_{i-1},\tau_{i}].
\end{eqnarray}
To see all these Ricci flows are well-concatenated, we apply (\ref{breather_condition}) to observe that
\begin{eqnarray*}
g_1(\tau_0)&=&\alpha^{-1}\phi^*g_0(0)=g_0(1),
\\
g_i(\tau_{i-1})&=&\alpha^{-i}(\phi^i)^*g_0(0)=\alpha^{-(i-1)}(\phi^{i-1})^*g_0(1)
\\
&=&\alpha^{-(i-1)}(\phi^{i-1})^*g_0\left(\alpha^{i-1}(\tau_{i-1}-\tau_{i-2})\right)=g_{i-1}(\tau_{i-1}).
\end{eqnarray*}
Therefore we define an ancient solution
\begin{eqnarray}\label{ancient_solution}
g(\tau)=\left\{ \begin{array}{ll}
g_0(\tau) & \text{for } \tau\in[0,1] \\
g_i(\tau) & \text{for } \tau\in[\tau_{i-1},\tau_{i}]\text{ and  } i\geq 1.\end{array}\right.
\end{eqnarray}
It then follows from the uniqueness theorem of Chen and Zhu \cite{chen2006uniqueness} that the ancient solution $g(\tau)$ is smooth.
\\

Now we proceed to show that $(M,g(\tau))_{\tau\in[0,\infty)}$, where $g(\tau)$ is defined in (\ref{ancient_solution}), is Type I. We need only to consider the case when $\tau\geq 1$. Let $i\geq 1$ be such that $\tau\in[\tau_{i-1},\tau_i]$. Then
\begin{eqnarray*}
|Rm(g(\tau))|=|Rm(g_i(\tau))|\leq\alpha^{i}\sup_{M\times[0,1]}\left|Rm\Big((\phi^i)^*g_0(\tau)\Big)\right|\leq C\alpha^{i},
\end{eqnarray*}
where we have used (\ref{bounded_condition}), (\ref{piece_metric}), and (\ref{ancient_solution}). Then we have
\begin{eqnarray} \label{Type_I_2}
 |Rm(g(\tau))|\leq C\alpha^{i}\leq\frac{C}{\tau}\tau_i\alpha^{i}\leq \frac{B}{\tau},
\end{eqnarray}
where we have used (\ref{time_bound}), and $B=CC_0$ is independent of $i$.
\\

With all the preparations, we are ready to prove our main theorem.

\begin{proof}[Proof of Theorem \ref{Theorem_Main}]
Fix an arbitrary point $y\in M$ as the base point, and for each $i\geq 0$ define
\begin{eqnarray}\label{x_i_def}
 x_i=\phi^{-(i+1)}(y).
\end{eqnarray}
In Lu and Zheng \cite{lu2017new}, they made an assumption that $\{x_i\}_{i=1}^\infty$ are not drifted away to space infinity so as to apply Theorem 4.1 in \cite{cao2011conjugate} to show that $\{(M,\tau_i^{-1}g(\tau\tau_i),(x_i,1))_{\tau=[1,2]}\}_{i=1}^\infty$ converges, after passing to a subsequence, to the canonical form of a shrinking gradient Ricci soliton. Instead we will show that $l(x_i,\tau_i)$, where $i\geq 0$ and $l$ is the reduced distance based at $(y,0)$,  is a bounded sequence. To see this, we let $\sigma:[0,1]\rightarrow M$ be a smooth curve such that $\sigma(0)=y$ and $\sigma(1)= x_0$. Let $A<\infty$ be such that
\begin{eqnarray}\label{velocity}
|\dot{\sigma}(\tau)|_{g_0(\tau)}\leq A, \text{ for all } \tau\in[0,1].
\end{eqnarray}
For each $i\geq 0$, we define
\begin{eqnarray}\label{sigmadef}
\sigma_i(\tau):=\phi^{-(i+1)}\circ\sigma(\alpha^{i+1}(\tau-\tau_{i})),\ \text{where } \tau\in[\tau_i,\tau_{i+1}].
\end{eqnarray}
We observe that these $\sigma_i$'s and $\sigma$ altogether define a piecewise smooth curve in $M$:
\begin{eqnarray*}
\sigma_0(\tau_0)&=&\phi^{-1}\circ\sigma(0)=\phi^{-1}(y)=x_0=\sigma(1),
\\
\sigma_i(\tau_i)&=&\phi^{-(i+1)}\circ\sigma(0)=\phi^{-i}\circ\sigma(1)
\\
&=&\phi^{-i}\circ\sigma(\alpha^i(\tau_i-\tau_{i-1}))=\sigma_{i-1}(\tau_i).
\end{eqnarray*}
We then define $\gamma_i:[0,\tau_{i+1}]\rightarrow M$, where $i\geq 0$, as
\begin{eqnarray*}
\gamma_i(\tau):=\left\{\begin{array}{ll}
\sigma(\tau) & \text{when } \tau\in[0,1],
\\
\sigma_j(\tau) & \text{when } \tau\in[\tau_j,\tau_{j+1}] \text{ and } 0\leq j\leq i.
\end{array}\right.
\end{eqnarray*}
Apparently $\gamma_i(\tau)$ is piecewise smooth, and $\gamma_i(0)=y$, $\gamma_i(\tau_{i+1})=\phi^{-(i+2)}(y)=x_{i+1}$. We compute for $i\geq 0$
\begin{eqnarray*}
\mathcal{L}(\gamma_i)&=&\mathcal{L}(\sigma)+\sum_{j=0}^i\int_{\tau_j}^{\tau_{j+1}}\sqrt{\tau}\Big( R(\sigma_j(\tau),\tau)+|\dot{\sigma_j}(\tau)|^2_{g_{j+1}(\tau)}\Big)d\tau
\\
&\leq& D + \sum_{j=0}^i\int_{\tau_j}^{\tau_{j+1}}\sqrt{\tau}\Big( \frac{B}{\tau}+A\alpha^{j+1}\Big)d\tau
\end{eqnarray*}
where in the last inequality we have used $D$, a constant independent of $i$, to represent $\mathcal{L}(\sigma)$, and we have used the Type I condition (\ref{Type_I_2}), the definition  (\ref{sigmadef}) of $\sigma_j$, and the assumption (\ref{velocity}). Continuing the computation using (\ref{time_bound}), we have
\begin{eqnarray*}
\mathcal{L}(\gamma_i)\leq D+C_1\sum_{j=0}^i \alpha^{-\frac{j+1}{2}},
\end{eqnarray*}
where $C_1$ is a constant independent of $i$. It follows from the definition (\ref{reduced_distance}) that
\begin{eqnarray*}
l(x_{i+1},\tau_{i+1})&\leq&\frac{1}{2\sqrt{\tau_{i+1}}}\mathcal{L}(\gamma_i)
\\
&\leq&\frac{1}{2}D\alpha^{\frac{i+1}{2}}+\frac{1}{2}C_1\sum_{j=0}^i\alpha^{\frac{j}{2}}\leq C_2<\infty,
\end{eqnarray*}
where $C_2$ is a constant independent of $i$, and we have used $\alpha^{\frac{1}{2}}\in(0,1)$.
\\

Now we consider the sequence
\begin{eqnarray}\label{sequence}
\{(M,\tau_{i}^{-1}g(\tau\tau_i),(x_i,1))_{\tau\in[1,\alpha^{-1}]}\}_{i=1}^\infty.
\end{eqnarray}
We observe that
\begin{eqnarray*}
\tau_{i}^{-1}g(\tau_i)=\tau_{i}^{-1}\alpha^{-(i+1)}\Big(\phi^{i+1}\Big)^*g_0(0),
\end{eqnarray*}
where $\tau_{i}^{-1}\alpha^{-(i+1)}$ is bounded from above and below by constants independent of $i$, because of  (\ref{time_bound}). Taking into account the definition (\ref{x_i_def}) of $x_i$, we can use
\begin{eqnarray*}
\operatorname{inj}_{g_0(0)}(y)>0
\end{eqnarray*}
to verify the condition (\ref{inj}). It follows from Theorem \ref{asymoptotic_shrinker} that (\ref{sequence}) converges smoothly to the canonical form of a shrinking gradient Ricci soliton. Furthermore, since $(M,\tau_{i}^{-1}g(\tau_i),x_i)$ and $(M,g_0(0), y)$ differ only by a bounded scaling constant and a diffeomorphism that preserves the base points, by the definition of the Cheeger-Gromov convergence, such diffeomorphism does not affect the limit. In other words, there exists a constant $C_3>0$, such that
\begin{eqnarray*}
(M,\tau_{i}^{-1}g(\tau_i),x_i)\rightarrow (M,C_3g_0(0),y)
\end{eqnarray*}
in the pointed smooth Cheeger-Gromov sense. Therefore $(M,g_0(0),y)$ also has a shrinker structure up to scaling. It then follows from the backward uniqueness of Kotschwar \cite{kotschwar2010backwards} that the shrinking breather $(M,g_0(\tau))_{\tau\in[0,1]}$ is the canonical form of a shrinking gradient Ricci soliton.

\end{proof}

\textbf{Acknowledgement:} The author is grateful to Professor Peng Lu and Professor Qi Zhang for their interest in this problem.

\bibliographystyle{plain}
\bibliography{citation}

\begin{thebibliography}{1}

\bibitem{cao2011conjugate}
Xiaodong Cao and Qi~S. Zhang.
\newblock The conjugate heat equation and ancient solutions of the {Ricci}
  flow.
\newblock {\em Advances in Mathematics}, 228(5):2891--2919, 2011.

\bibitem{chen2006uniqueness}
Bing-Long Chen and Xi-Ping Zhu.
\newblock Uniqueness of the {Ricci} flow on complete noncompact manifolds.
\newblock {\em Journal of Differential Geometry}, 74(1):119--154, 2006.

\bibitem{chow2007ricci}
Bennett Chow, Sun-Chin Chu, David Glickenstein, Christine Guenther, Jim
  Isenberg, Tom Ivey, Dan Knopf, Peng Lu, Feng Luo, and Lei Ni.
\newblock {\em The {Ricci} {Flow}: {Techniques} and {Applications}: {Part} {I}:
  {Geometric} {Aspects}}.
\newblock American Mathematical Society, 2007.

\bibitem{hamilton1995compactness}
Richard Hamilton.
\newblock A compactness property for solutions of the {Ricci} flow.
\newblock {\em American Journal of Mathematics}, 117(3):545--572, 1995.

\bibitem{kotschwar2010backwards}
Brett Kotschwar.
\newblock Backwards uniqueness for the ricci flow.
\newblock {\em International Mathematics Research Notices},
  2010(21):4064--4097, 2010.

\bibitem{lu2017new}
Peng Lu and Yu~Zheng.
\newblock New proofs of {Perelman's} theorem on shrinking breathers in {Ricci}
  flow.
\newblock {\em The Journal of Geometric Analysis}, pages 1--7, 2017.

\bibitem{naber2010noncompact}
Aaron Naber.
\newblock Noncompact shrinking four solitons with nonnegative curvature.
\newblock {\em Journal f{\"u}r die reine und angewandte Mathematik (Crelles
  Journal)}, 2010(645):125--153, 2010.

\bibitem{perelman2002entropy}
Grisha Perelman.
\newblock The entropy formula for the {Ricci} flow and its geometric
  applications.
\newblock {\em arXiv preprint math/0211159}, 2002.

\bibitem{zhang2014no}
Qi~S Zhang.
\newblock A no breathers theorem for some noncompact {Ricci} flows.
\newblock {\em Asian Journal of Mathematics}, 18(4):727--756, 2014.

\end{thebibliography}

\noindent Department of Mathematics, University of California, San Diego, CA, 92093
\\ E-mail address: \verb"yoz020@ucsd.edu"

\end{document}